\documentclass[a4paper,12pt] {article}
\usepackage{etex}

\usepackage{graphics}
\usepackage{graphicx}
\usepackage{epsfig}
\usepackage{pst-grad} 
\usepackage{pst-plot} 
\usepackage{pstricks}
\usepackage{lmodern}
\usepackage{xcolor,rotating,epic,eepic}
\usepackage[permil]{overpic}
\usepackage[english]{babel}
\usepackage{lscape}
\usepackage{array}
\usepackage{multirow}
\usepackage{mathrsfs}
\usepackage{amsmath}
\usepackage{amssymb}
\usepackage{float}
\usepackage{subfig}
\usepackage{multirow}
\usepackage{placeins}
\usepackage{cite}

\definecolor{lgray}{gray}{0.85}
\definecolor{myorange}{rgb}{0.87,0.49,0.0}

\begin{document}
\begin{center}
\Large{\textbf{Nonlinear Normal Modes, Modal Interactions and Isolated Resonance Curves}}\vspace{0.5cm}
\end{center}

\begin{center}
\large{\textbf{R.J. Kuether$^1$, L. Renson$^2$, T. Detroux$^2$, C. Grappasonni$^2$, \\ G. Kerschen$^2$, M.S. Allen$^1$}}\\\vspace{0.5cm}
\normalsize
$^1$ Department of Engineering Physics, \\
University of Wisconsin-Madison,\\
Madison, WI 53706.\\\vspace{3mm}
$^2$Space Structures and Systems Laboratory,\\
Department of Aerospace and Mechanical Engineering,\\
University of Li\`ege, Li\`ege, Belgium. \\
\vspace{0.35cm} {Corresponding author}: L. Renson\\
Email: l.renson@ulg.ac.be, phone: +32 4 3664854.
\vspace{0.35cm}

\begin{abstract}
The objective of the present study is to explore the connection between the nonlinear normal modes of an undamped and unforced nonlinear system and the isolated resonance curves that may appear in the damped response of the forced system. To this end, an energy balancing technique is used to predict the amplitude of the harmonic forcing that is necessary to excite a specific nonlinear normal mode. A cantilever beam with a nonlinear spring at its tip serves to illustrate the developments. The practical implications of isolated resonance curves are also discussed by computing the beam response to sine sweep excitations of increasing amplitudes.\\
\vspace{0.35cm}
\end{abstract}

\end{center}
\normalsize

\section{Introduction}\label{sec_intro}

Nonlinearity is important in many structural dynamic applications that are of interest to engineers, for example in structures with bolted interfaces \cite{Segalman}, machinery with rubber isolation mounts, microsystems subjected to thermal, magnetic or friction forces \cite{Czaplewski}, and biomechanical systems \cite{Garcia}.  In other cases the baseline structure is linear, but its performance can be enhanced by adding or engineering certain types of nonlinearities \cite{Gendelman,VakakisJVA}.  However, nonlinear dynamics is a rich and complicated field and so engineers tend instead to ignore nonlinearity or to seek a linear model that approximates the system at the forcing level of interest \cite{Segalman2}.

Vibration modes form the foundation of our understanding of linear dynamic systems, and influence efforts related to testing, modeling, validating and controller design.  Rosenberg \cite{Rosenberg} extended modal analysis to nonlinear systems in the 1960's, coining the phrase nonlinear normal mode (NNM).  The area received new attention in the 1990's \cite{Shaw,VakakisMSSP,VakakisBook} and now it is clear that NNMs can be used to obtain a wealth of insight into the response of a nonlinear system \cite{VakakisMSSP,Kerschen}.  For example, NNMs have been used to explain internally resonant and non-resonant motions of structures \cite{Lacarbonara}, to design a nonlinear vibration absorber (also called a nonlinear energy sink) \cite{Gendelman}, to create or validate a reduced order model for a system \cite{Touze}, and to explain changes in the oscillation frequency and the deformation shape of the free and forced response of a structure \cite{Kerschen}.
	
In recent years important progress has been made in the numerical calculation of undamped \cite{Arquier,Kuether,Laxalde,Peeters} and damped NNMs \cite{Pesheck,Renson}.  These new algorithms have been used to compute the nonlinear modes of a geometrically nonlinear finite element model of a component from a diesel exhaust system, a full-scale aircraft, a bladed disk from a turbine and a strongly nonlinear satellite in \cite{Kuether,Kerschen2,Krack,Renson2}, respectively.  A framework for experimental identification of NNMs was recently presented in \cite{Peeters2} and validated on an academic structure \cite{Peeters3}.  More recent works have begun to use this framework on more complicated structures \cite{Allen,Zapico}.
	
One fundamental property of undamped NNMs is the fact that they can be realized when a harmonic forcing function cancels the damping force in the damped system \cite{Peeters2}.  As a result, they form the backbone of the nonlinear forced response (NLFR) curves \cite{VakakisBook,Kerschen,Wagg} and hence they approximate the oscillation frequency and deformation shape that are exhibited at resonance, when a structure is at the greatest risk of failure.  The relationship between the NLFR and the NNM backbone is simple for mild nonlinearities, but most realistic systems exhibit complicated NNMs with many interactions between the various modes leading to internally resonant branches. Many works have explored interactions between nonlinear modes with commensurate linear frequencies, e.g., when a pair of the linearized natural frequencies of the system have an integer ratio \cite{NayfehB,Wagg2,Rega}. Some of them have even exploited these modal interactions for optimal design \cite{Dou}. In contrast, only a few have considered the case that is of interest in this work where the modal interactions occur between pairs of modes whose linear frequencies are not integer related. The investigation of such interactions requires one to resort to computational methods, which is another specific aspect of this study.

This work explores the relationships between these interacting NNMs and the forced response of the nonlinear system, especially for the case in which the forced response shows an isolated resonance curve (IRC). Specifically, we show that the interactions between NNMs are responsible for the IRCs in the forced response. We note that the paper \cite{Neild} discusses the relationships between bifurcations of backbone curves and IRCs; it is therefore the ideal companion of the present study. These detached families of solutions are frequently not detected because they do not emerge naturally from the fundamental response when numerical continuation is used. They can lie outside or inside the main resonance curve \cite{Gatti,Gatti2}, with the former case typically being more important because one is likely to underestimate the response of the nonlinear system \cite{Alexander,Duan}. Isolated resonances may also go undetected during laboratory experiments when stepped/swept sine testing is employed.

The paper is organized as follows. Section \ref{Periodic} briefly reviews the methodology used to compute the periodic motions of the undamped and damped form of the nonlinear equations of motion, along with a phase resonance condition extended to nonlinear systems. An adaptation of the energy balancing procedure presented by Neild et al.~\cite{NeildE1,NeildE2}, which can be used to estimate the forcing amplitude required to isolate the NNM, is also presented in Section \ref{Periodic}. Section~\ref{NNMsec} applies the energy balancing technique to a nonlinear cantilever beam with the aim to predict the resonances from the knowledge of the NNMs and of the damping matrix. The predictions of  Section~\ref{NNMsec} are validated using the computation of NLFRs and bifurcation tracking in Section~\ref{NLFRsec}. Section~\ref{sec5} discusses the practical implications of IRCs by computing the beam response to sine sweep excitations of increasing amplitudes. The conclusions of the present study are presented in Section~\ref{sec_future}.

\section{Periodic Solutions of a Nonlinear System}\label{Periodic}

\subsection{Computation of Nonlinear Normal Modes and Nonlinear Frequency Responses}\label{Sec21}

The $N$-degree-of-freedom (DOF) equations of motion for a nonlinear system generally can be written as
\begin{equation}\label{EOM}
  \mathbf{M\ddot{x}}+\mathbf{C\dot{x}}+\mathbf{Kx}+\mathbf{f}_{NL}(\mathbf{x})=\mathbf{f}(t)
\end{equation}
	 	
The $N\times N$ matrices $\mathbf{M}$, $\mathbf{C}$, and $\mathbf{K}$ represent the mass, damping and stiffness matrices, respectively. The displacement, velocity and acceleration are represented with the $N\times 1$ vectors $\mathbf{x}$, $\mathbf{\dot{x}}$, and $\mathbf{\ddot{x}}$, and the external loads are applied through the $N\times 1$ force vector $\mathbf{f}(t)$. The $N\times 1$ nonlinear restoring force vector, $\mathbf{f}_{NL}(\mathbf{x})$, accounts for the nonlinearities in the physical system. We only consider the case where the nonlinear restoring force depends on displacement.

The NNMs calculated in this study are defined as {\it not necessarily synchronous periodic motions of the undamped and unforced nonlinear system} \cite{Kerschen,Lee}. A variety of methods exist to find these periodic solutions, e.g., perturbation techniques \cite{Gendelman} and harmonic balance \cite{Cochelin}. In this paper, NNMs are computed using the shooting technique combined with numerical continuation \cite{Kuether,Peeters}. The outcome of the calculations is a frequency-energy plane which depicts the evolution of the fundamental frequency of the NNM as the energy changes.

Nonlinear forced responses (NLFRs), i.e., the periodic responses of the damped system to a monoharmonic excitation force
\begin{equation}\label{force}
  \mathbf{f}(t)=\mathbf{A} \sin \left( \omega t\right)
\end{equation}
are also calculated herein using shooting and numerical continuation~\cite{Sracic}. NLFRs can reveal new phenomena that cannot be observed with linear theory, such as frequency-energy dependence, subharmonic and superharmonic resonances, coexisting solutions, and stable/unstable periodic motions.

\subsection{Connection between Nonlinear Normal Modes and Nonlinear Frequency Responses}\label{Sec22}

\subsubsection{Nonlinear Phase Lag Quadrature Criterion}

The first theoretical connection between NLFRs and NNMs can be derived thanks to the nonlinear phase lag quadrature criterion \cite{Peeters2}. Specifically, the damped system can be made to respond in a single NNM motion if the excitation exactly cancels out the damping force. To this end, a multi-point, multi-harmonic excitation is to be applied to the system so that the harmonics of the response are all in quadrature with the harmonics of the excitation.

In practice, it is unlikely that the forcing of interest will exactly cancel damping as outlined above, because this requires that the force be distributed in space and that it be comprised of many harmonics. A few studies \cite{Peeters2,Peeters3,Allen,Zapico,Kuether2} have shown that an accurate approximation to the NNM can be obtained by a much simpler force that excites resonance. For example, Peeters et al. \cite{Peeters2,Peeters3} explored whether a single-point, monoharmonic excitation could approximately isolate an NNM, and found good results in simulations and experiment with a lightly damped beam. In their efforts it was helpful to define a multi-harmonic mode indicator function (MIF) which indicates when the 90 degrees phase lag condition has been obtained.  When a single-point sinusoidal force is applied to the nonlinear structure, which is the case considered in this paper, the MIF is defined as follows
\begin{equation}\label{MIF}
  \Delta_1=\frac{\mbox{Re}(\mathbf{Z}_1)^*\mbox{Re}(\mathbf{Z}_1)}{\mathbf{Z}_1^*\mathbf{Z}_1}
\end{equation}
where the operator $(\cdot)^*$ represents the complex conjugate transpose and $\mathbf{Z}_1$ is the complex Fourier coefficient of the fundamental harmonic of the computed NLFR. The MIF in Eq. (\ref{MIF}) indicates that resonance occurs when $\Delta_1$ is equal to one.

\subsubsection{An Energy Balancing Technique}\label{Sec222}
The relationship between NLFRs and NNMs was also studied using an energy balancing technique in \cite{NeildE1,NeildE2}. Based on the second-order normal form theory, the analytical developments hold for weakly nonlinear regimes of motion. This technique is slightly revisited herein by employing a numerical viewpoint, which allows one to consider more strongly nonlinear regimes.

Let us first consider a linear system. As shown in \cite{Geradin}, if the system is oscillating in a linear normal mode denoted as $\mathbf{x}(t)$, then the damping forces instantaneously exert a distributed force $\mathbf{C\dot{x}(t)}$ and the power dissipated at any instant is
\begin{equation}\label{power}
  P_{diss}=\mathbf{\dot{x}}(t)^{\mbox{T}}\mathbf{C\dot{x}(t)}
\end{equation}
and the total energy dissipated over one cycle is
\begin{equation}\label{energy}
  E_{diss/cyc}=\int_0^T P_{diss}\,dt
\end{equation}
Similarly, an arbitrary forcing function $\mathbf{f}(t)$ inputs energy into the system as
\begin{equation}\label{energy2}
  E_{in/cyc}=\int_0^T \mathbf{\dot{x}}(t)^{\mbox{T}} \mathbf{f}(t)\,dt
\end{equation}
At resonance, the energy dissipated by the damping forces must match the total energy input to the system over the period $T$. The balance is enforced by setting $E_{diss/cyc}=E_{in/cyc}$ \cite{Geradin}. For a single-point, monoharmonic force with amplitude $A$, the scaling on $A$ can be computed by satisfying
\begin{equation}\label{puissancereactive}
  A\int_0^T \mathbf{\dot{x}}(t)^{\mbox{T}}\left(\mathbf{e}_n \sin{i\omega t}\right)\,dt=\int_0^T \mathbf{\dot{x}}(t)^{\mbox{T}}\mathbf{C\dot{x}(t)} \,dt
\end{equation}
where $\mathbf{e}_n$ is a vector of zeros with a value of one at the location $n$, which is the point at which the force is applied. This energy balancing criterion is a useful result, because it enables the practitioner to establish formally the direct link from the computed linear normal modes, i.e., the periodic motions of the undamped, unforced system, to the resonant response of the damped forced system.

The energy balance, $E_{diss/cyc}=E_{in/cyc}$, also holds for nonlinear systems. So, if both the NNMs $\mathbf{x}(t)$ and the damping $\mathbf{C}$ in the system are known, Eq. (\ref{puissancereactive}) can be readily used to estimate the forcing amplitude $A$ that would excite the system at resonance with associated motion $\mathbf{x}(t)$. While it is common practice to excite a system using a monoharmonic force, one should note that higher harmonics might be necessary to achieve a reasonable approximation to the NNM motion, especially near internal resonances, so any calculations based on Eq. (\ref{puissancereactive}) should be regarded as approximate.

\section{Prediction of the Forced Response of a Cantilever Beam using Nonlinear Normal Modes}\label{NNMsec}

In the present and next sections, a model of a cantilevered beam with a cubic nonlinear spring attached at the beam tip is used to investigate the connection between NNMs and NLFRs.

\begin{figure}[t]
\begin{center}
\includegraphics[width=0.5\textwidth]{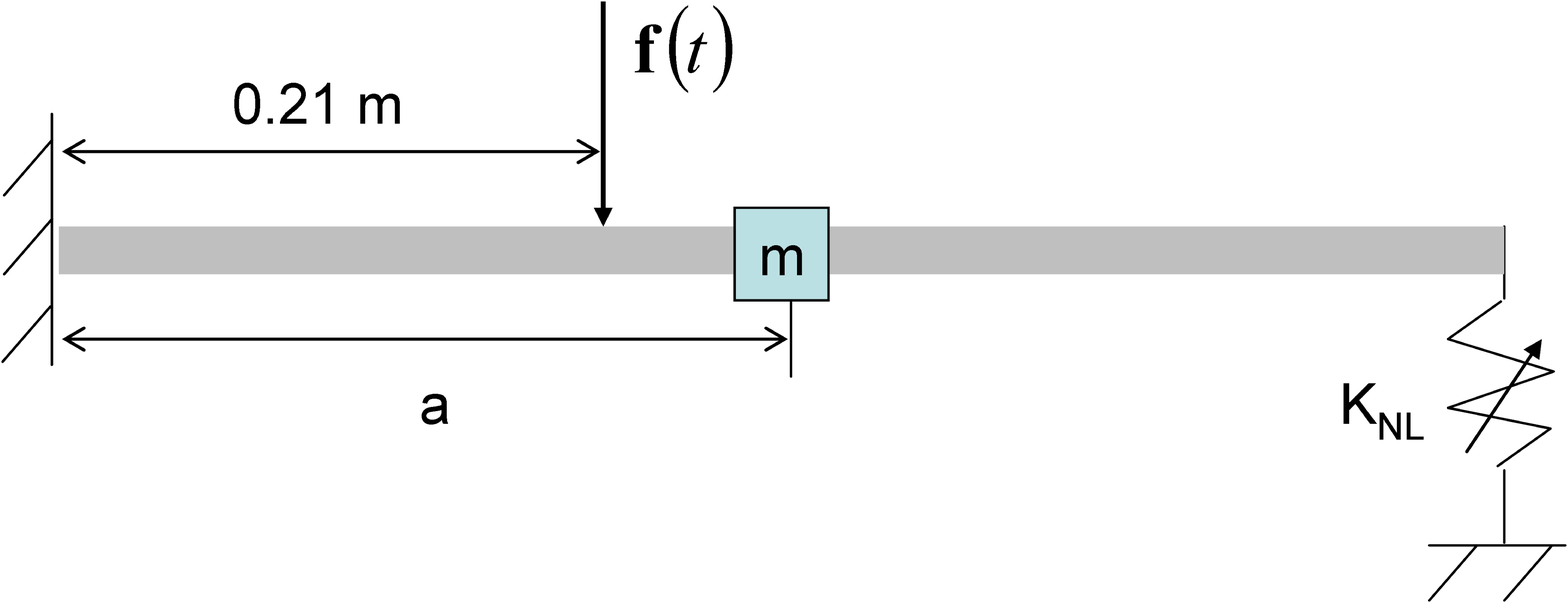}
\caption{Schematic of a cantilever beam with a cubic nonlinear spring attached to the beam tip and a modifying lumped mass of 0.5 kg. The addition of the mass shifted the location of the 3:1 modal interaction with NNM 2 on the first NNM branch (as seen later in Fig. \ref{fig2}(a)).}
\label{fig1}
\end{center}
\end{figure}

The beam was 0.7 m in length, with a width and thickness of 0.014 m, and was constructed of structural steel with a Young's modulus of 205 GPa and a density of 7800 kg/m$^3$. A schematic of the FEA model is shown in Fig. \ref{fig1}. A linear finite element model of the planar beam was created in Abaqus using 20 B31 Euler-Bernoulli beam elements, giving it a total of 60 DOF. A mass and stiffness proportional damping model was used, defining the damping matrix as $\mathbf{C}=a\mathbf{K}+b\mathbf{M}$ with $a=-0.0391$ and $b=1.47\,10^{-4}$. These parameters were chosen such that the damping ratios of the first and second linear modes were 1\% and 5\%, respectively. The cubic nonlinear spring had a coefficient of $K_{NL} = 6\,10^9 N/m^3$, and was attached at the beam tip affecting only the transverse direction.
\begin{figure}[!htp]
\begin{center}
\begin{tabular}{c}
\includegraphics[width=0.85\textwidth]{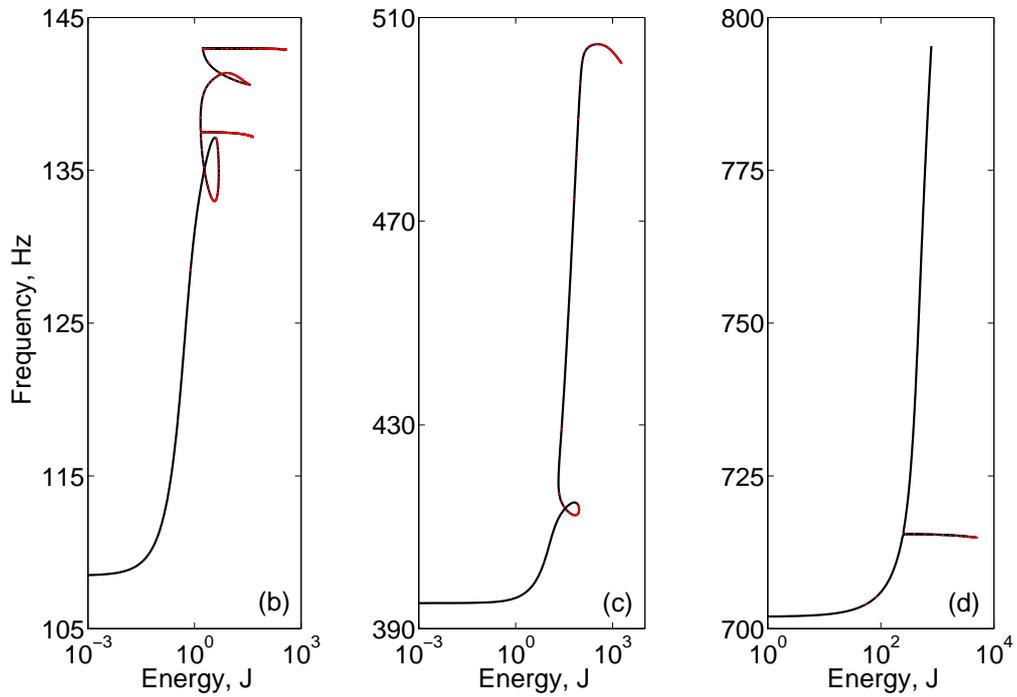}
\end{tabular}
\caption{The nonlinear normal modes of the nonlinear beam: (a) NNM 1, (b) NNM 2, (c) NNM 3, and (d) NNM 4. Solid black represents a stable solution, and dotted red represents an unstable solution. The detailed view of NNM 1 shows the crossing of higher order NNMs in dashed lines at fractions of their fundamental frequency.}
\label{fig2}
\end{center}
\end{figure}

The frequency-energy plots (FEPs) of the first four NNMs are shown in Fig. \ref{fig2}. These FEPs have two distinct features, namely a backbone, and tongues that emanate from the backbone. The backbone of NNM 1 in Fig. \ref{fig2}(a) shows an increase in fundamental frequency as the energy in the periodic solutions increased revealing that the nonlinear spring has a stiffening effect on this mode. Each of the first four NNMs showed this stiffening behavior, but the energy level at which the nonlinearity began to affect the frequency and deformation varied for each. The tongues that emerged from the backbones along each NNM are referred as modal interactions, or internal resonances, and occur when two or more NNMs interact. The response of the NNM at locations on these tongues showed a strong, multi-harmonic response at an integer frequency ratio of the interacting nonlinear mode.

Figure \ref{fig2}(a) displays a detailed view of the FEP of NNM 1, where the dashed, colored lines represent the frequency-energy behavior of the higher order NNMs after dividing the frequency by various integers. By shifting these NNMs down the frequency axis, it was possible to observe the location where the backbones of higher modes intersect with the NNM 1 backbone and cause a modal interaction to occur. Considering the modal interaction at approximately 37 Hz, which has the appearance of the Greek letter $\alpha$ and will hereafter be referred to as an $\alpha$-tongue, the 1/3rd frequency branch of the NNM 2 branch intersects the backbone of NNM 1. This causes NNM 1 to bifurcate and create a 3:1 internal resonance tongue that has solutions where NNM 1 and 2 interact. The other three modal interactions along NNM 1 were a 9:1 interaction with NNM 3 near 44 Hz, a 15:1 interaction with NNM 4 near 47 Hz, and a 13:1 interaction with NNM 4 near 54 Hz. It is important to note that NNM 1 in Fig. \ref{fig2}(a) was almost certainly incomplete, because in reality many more tongues could emanate from the backbone as the frequency is equal to many other integer fractions of higher NNMs. These additional tongues must have been missed by the continuation algorithm. In principle they could be found using a smaller stepsize but this becomes time consuming and was not pursued.

\begin{figure}[ht]
\begin{center}
\includegraphics[width=1.0\textwidth]{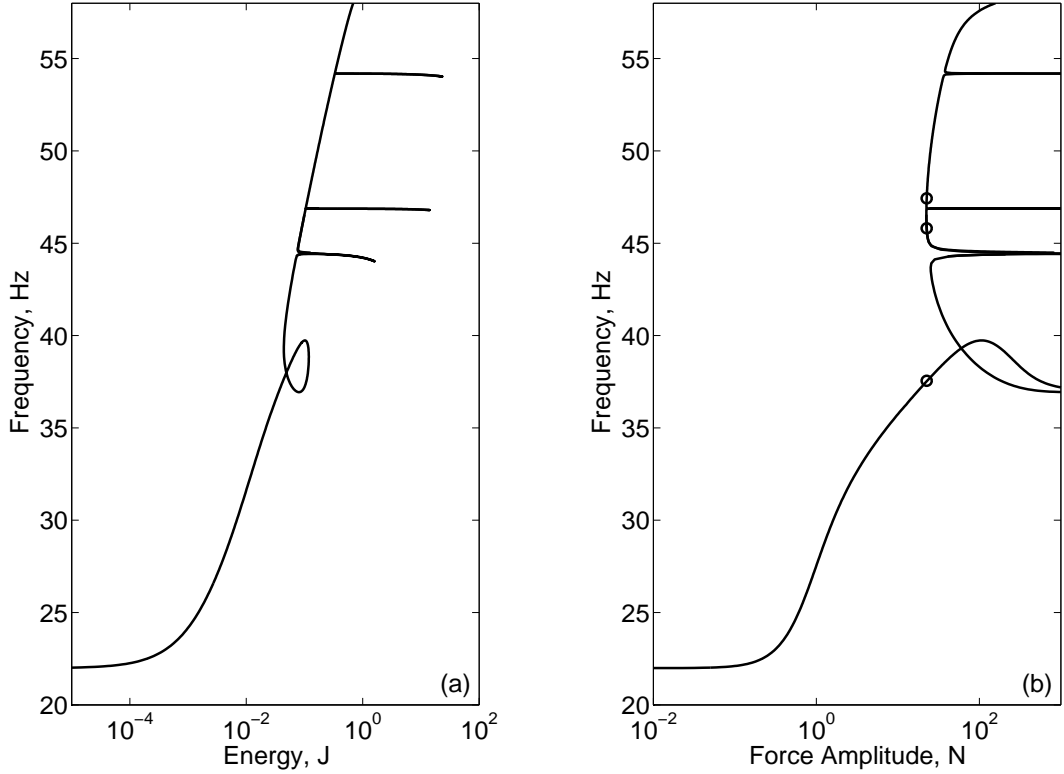}
\caption{First NNM: (a) FEP, (b) estimate of force amplitude required to obtain the motion given at each point on NNM 1. Circular markers indicate achievable resonance frequencies for a force of 22.6 N.}
\label{fig4}
\end{center}
\end{figure}

Focusing our attention on the first beam mode, our objective is to exploit the energy balancing technique of Section \ref{Sec222} to predict the system's nonlinear resonances based on the knowledge of the NNMs of Figure \ref{fig2}(a) and of the damping matrix. Specifically, Eq. (\ref{puissancereactive}) was used to estimate the monoharmonic driving force required to excite the NNM motion. The computed force amplitude and the corresponding frequency are displayed for the first NNM in Fig. \ref{fig4}(b); the NNM is repeated in Fig. \ref{fig4}(a). For forcing amplitudes smaller than 22.3 N, there is a unique solution meaning that, at a specific forcing amplitude, there exists a single resonance. For greater forcing amplitudes, this is no longer the case. For example, a force of 22.6 N could achieve multiple resonances at 37.5 Hz, 45.8 Hz and 47.4 Hz. One important remark is that the nonuniqueness of the resonances is due to the nonmonotonic increase in the forcing amplitude in Fig. \ref{fig4}(b), which is itself produced by the modal interactions (tongues) in the FEP in Fig. \ref{fig4}(a). 

\section{Calculation of the Forced Response of a Cantilever Beam using Numerical Continuation}\label{NLFRsec}

\subsection{Frequency Response Curves}\label{sec42}

\begin{figure}[ht]
\begin{center}
\includegraphics[width=0.75\textwidth]{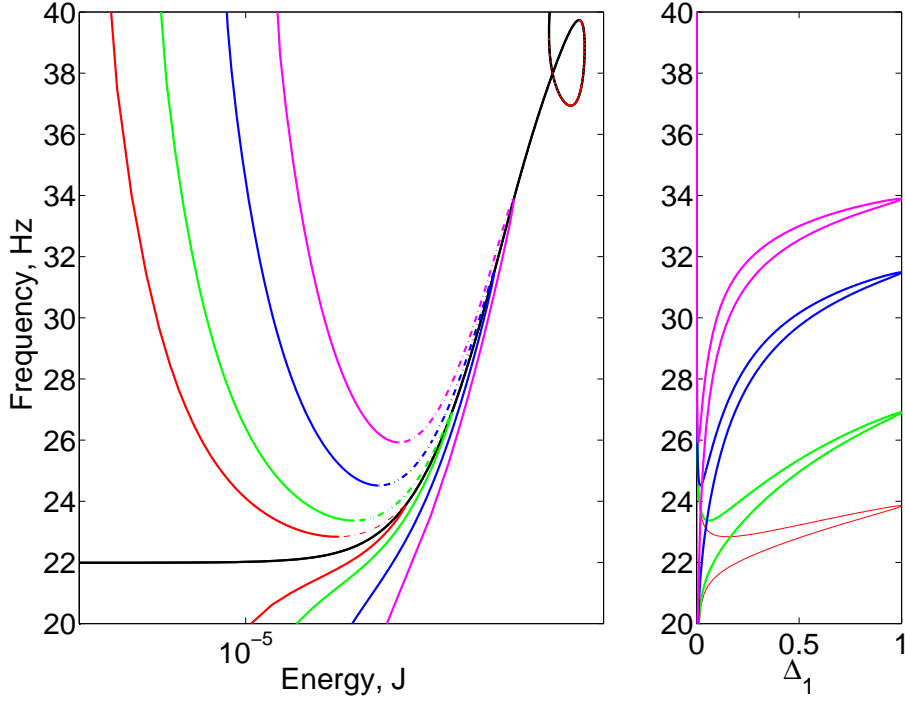}
\caption{(Left) NLFRs at frequencies near the first NNM where (solid) are stable periodic motions and (dash dot) are unstable periodic motions. The energy on the horizontal axis represents the maximum energy of each steady-state solution in the NLFR, and the vertical axis represents the forcing frequency. (Right) MIF of the forced response. The force amplitudes for each curve are (red) 0.445 N, (green) 0.890 N, (blue) 2.22 N, and (magenta) 4.45 N. The NNM is in black.}
\label{fig7}
\end{center}
\end{figure}

To validate the predictions made from Figure \ref{fig4}(b), the FEPs of the forced response were computed at different forcing amplitudes. Forcing amplitudes lower than the critical value of 22.3 N were first considered in Fig. \ref{fig7}. A classical behavior is observed in this figure where the forced response wraps around the NNM, acting as the backbone to the NLFR. The MIF from Eq. (\ref{MIF}) in the right plot in Fig. \ref{fig7} is approximately equal to 1 at resonance where a fold bifurcation changes the stability of the NLFR.

\begin{figure}
\begin{center}
\includegraphics[width=1.00\textwidth]{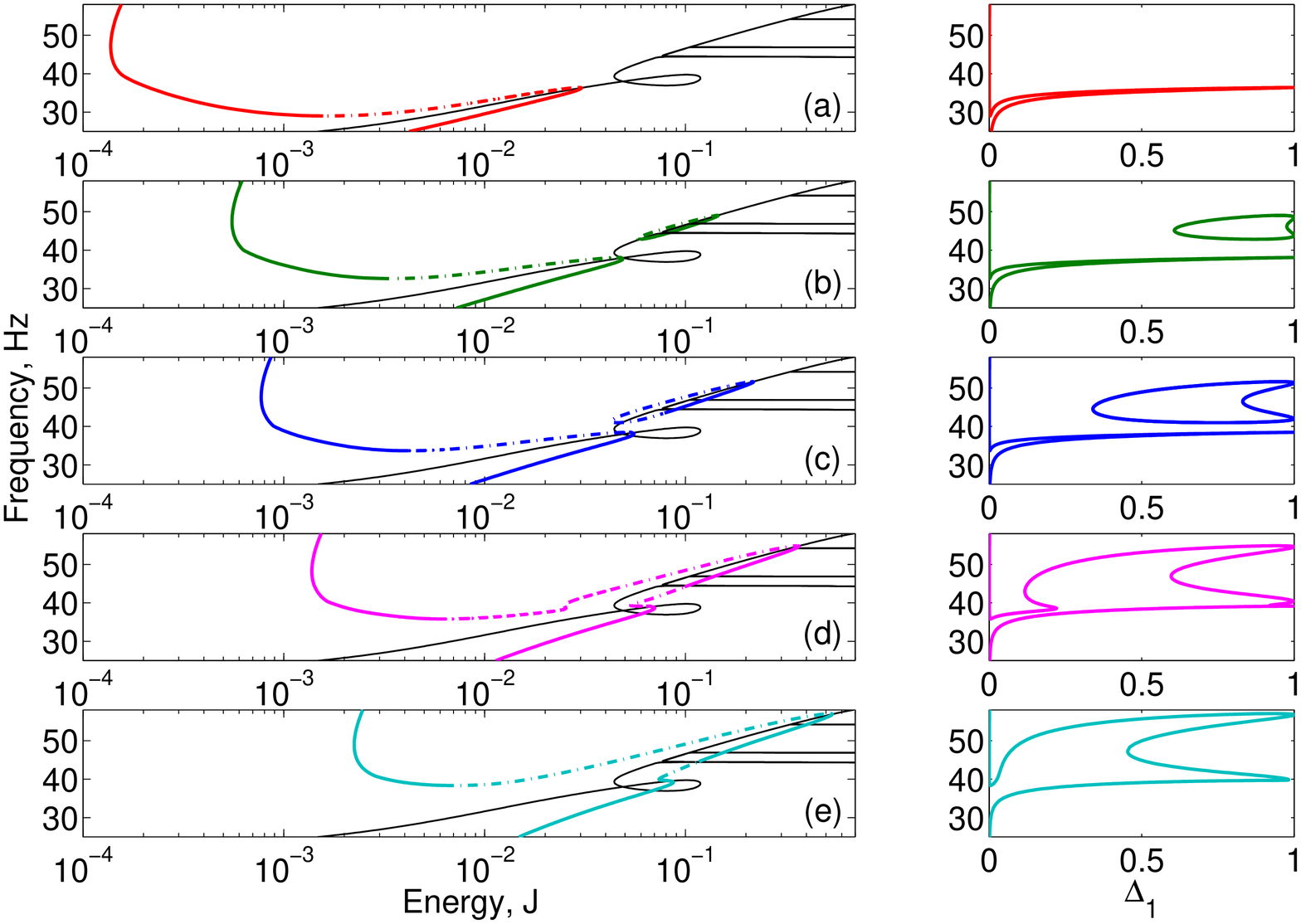}
\caption{(Left) NLFRs at frequencies near the first NNM where (solid) are stable periodic motions and (dash dot) are unstable periodic motions. (Right) MIF of the forced response. The force amplitudes for each curve (a - e) are (red) 11.3 N, (green) 22.6 N, (blue) 26.7 N, (magenta) 35.6 N, and (cyan) 45.2 N, respectively. . The NNM is in black.}
\label{fig8}
\end{center}
\end{figure}

Higher forcing amplitudes were considered in Fig. \ref{fig8}. For $A$ = 11.3 N in Fig. \ref{fig8}(a), the NLFR again wrapped around the backbone of the first NNM, as previously observed in Fig. \ref{fig7}. When the force amplitude doubled ($A$ = 22.6 N in Fig. \ref{fig8}(b)), three resonances were revealed by the MIF indicator at about 38.0 Hz (classical resonance), 44.2 Hz and 48.1 Hz. As clearly displayed in the NLFR, the two new resonances (associated with fold bifurcations) are responsible for the creation of an IRC. The response on this IRC is much larger than on the main branch so one would significantly underestimate the response if it was not detected. These results were found to be entirely consistent with the predictions of the energy balancing technique, which predicted multiple resonances from 22.3 N, and resonances at 37.5 Hz, 45.8 Hz and 47.4 Hz for 22.6 N. The fact that the IRC still wrapped around the backbone of NNM 1 underlines the essential role played by the NNM in the forced response. Fig. \ref{fig8}(b) also shows that the IRC was created when the NLFR branch approached the 3:1 interaction between NNM 1 and NNM 2; this observation suggests that interactions between NNMs are one possible driving mechanism for IRC onset.

\begin{figure}
\begin{center}
\includegraphics[width=1.00\textwidth]{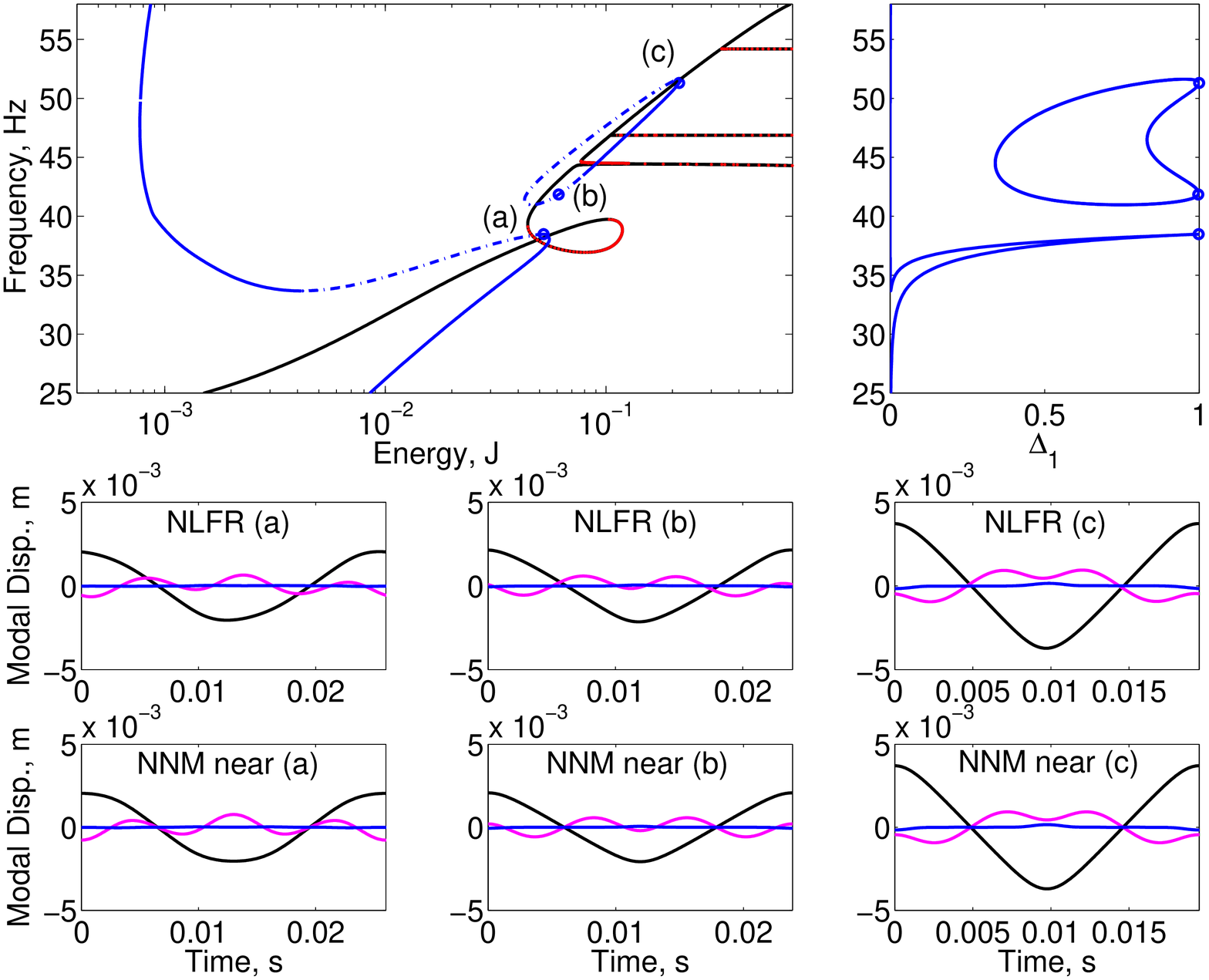}
\caption{Time histories projected onto unit displacement normalized modal coordinates at the three resonant conditions of the NLFR for $A$ = 26.7 N. The color code for the time response is (black) mode 1, (magenta) mode 2 and (blue) mode 3.}
\label{fig9}
\end{center}
\end{figure}

The forcing amplitude was slightly increased in Fig. \ref{fig8}(c) ($A$ = 26.7 N), and the resonant frequency on the main branch did not shift very much (from 38.1 Hz to 38.5 Hz). The IRC persisted and became larger, increasing its frequency range from 41 Hz to 51 Hz. A stable portion of the IRC has however become unstable through the emergence of Neimark-Sacker bifurcations (not represented). The latter will persist at higher forcing levels and are responsible for the quasiperiodic oscillations discussed in Section~\ref{sec5}. Further evidence of the connection between IRCs and modal interactions is given in Fig. \ref{fig9} which compares the time histories of the three resonant solutions (represented by circle markers) along the NLFR curve together with the time histories of the NNM near them. The periodic responses were projected onto the linear modal coordinates to better highlight the modal interaction. Indeed, the first two linear modal coordinates dominated all three of the resonant responses, and the comparison between the NLFR and the corresponding NNM all were in good agreement. There was slight phase shift however for solution (a), which can be explained by the use of a monoharmonic force input, whereas a multi-harmonic force would be needed to exactly isolate the NNM.

Increasing the amplitude even more ($A$ = 35.6 N in Fig. \ref{fig8}(d)) caused the main branch and the IRC to merge together, forming one continuous NLFR branch up to 55 Hz. The merging of these two branches therefore leads to a sudden and substantial change in the resonant frequency. For the highest forcing amplitude ($A$ = 45.2 N in Fig. \ref{fig8}(e)), the resonant frequency shifted to 56.6 Hz. This smooth increase in the resonant frequency continued at higher forces as well. Moving from Fig. \ref{fig8}(d) to Fig. \ref{fig8}(e) also caused one resonance to be eliminated, as indicated by the MIF.

\begin{figure}[tp]
\begin{center}
\includegraphics[width=0.75\textwidth]{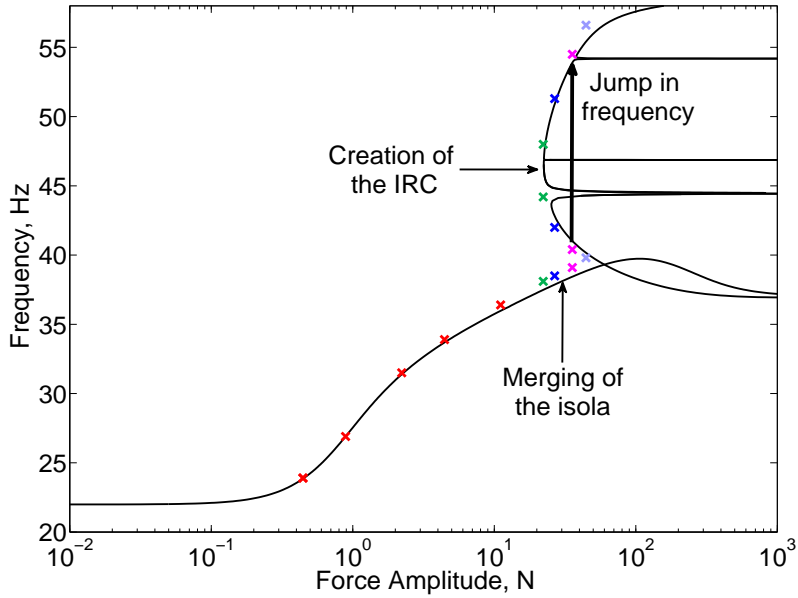}
\caption{Comparison between the predictions of the energy balancing criterion (solid black) and the forced resonance where the MIF is equal to one in Figs. \ref{fig7} and \ref{fig8} (cross markers).}
\label{fig10}
\end{center}
\end{figure}

For further validation of the energy balancing criterion (\ref{puissancereactive}), Figure \ref{fig10} superposes the responses of Figs. \ref{fig7} and \ref{fig8} where the MIF is equal to 1 onto Fig. \ref{fig4}(b). The cross markers, which represent the forced resonant response, are in close, though not exact, agreement with the predictions of the energy balancing criterion.

\subsection{Fold Bifurcations}

Because IRCs possess fold bifurcations, bifurcation tracking in the codimension-2 space (frequency-forcing amplitude-energy) is another tool that can reveal their existence. The procedure used in this section is based on the harmonic balance method described in \cite{Detroux}. 

\begin{figure}[!htp]
\begin{center}
\begin{tabular}{c}
\includegraphics[width=1.0\textwidth]{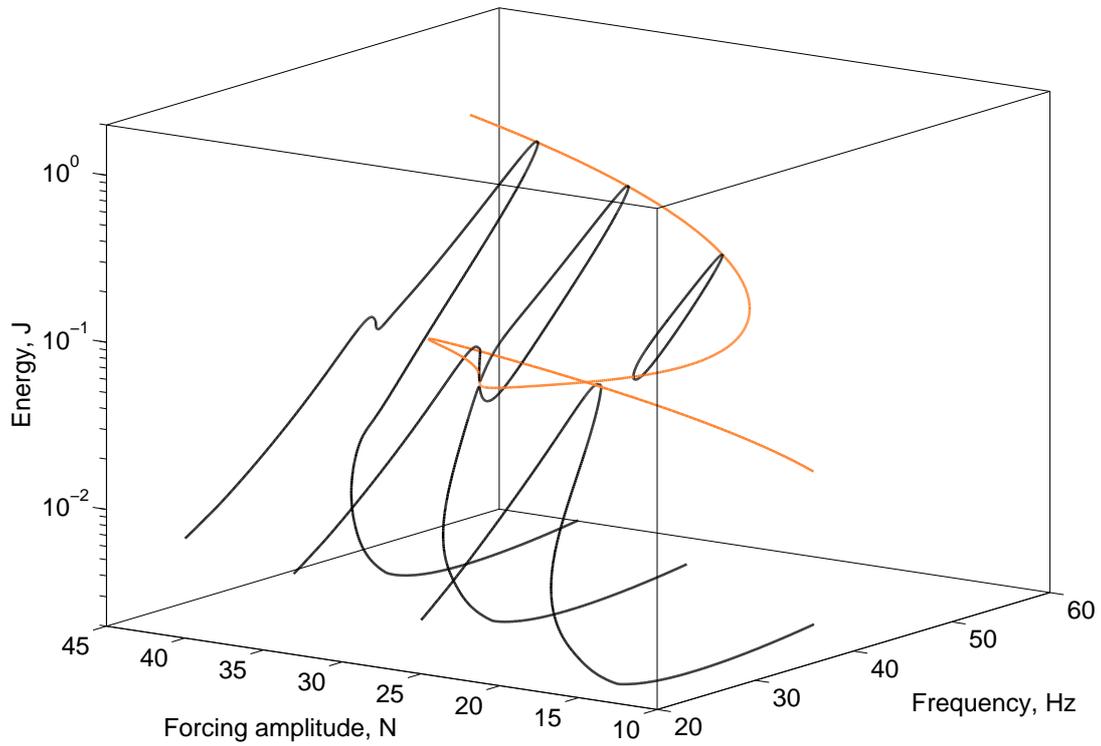}\\
\includegraphics[width=0.75\textwidth]{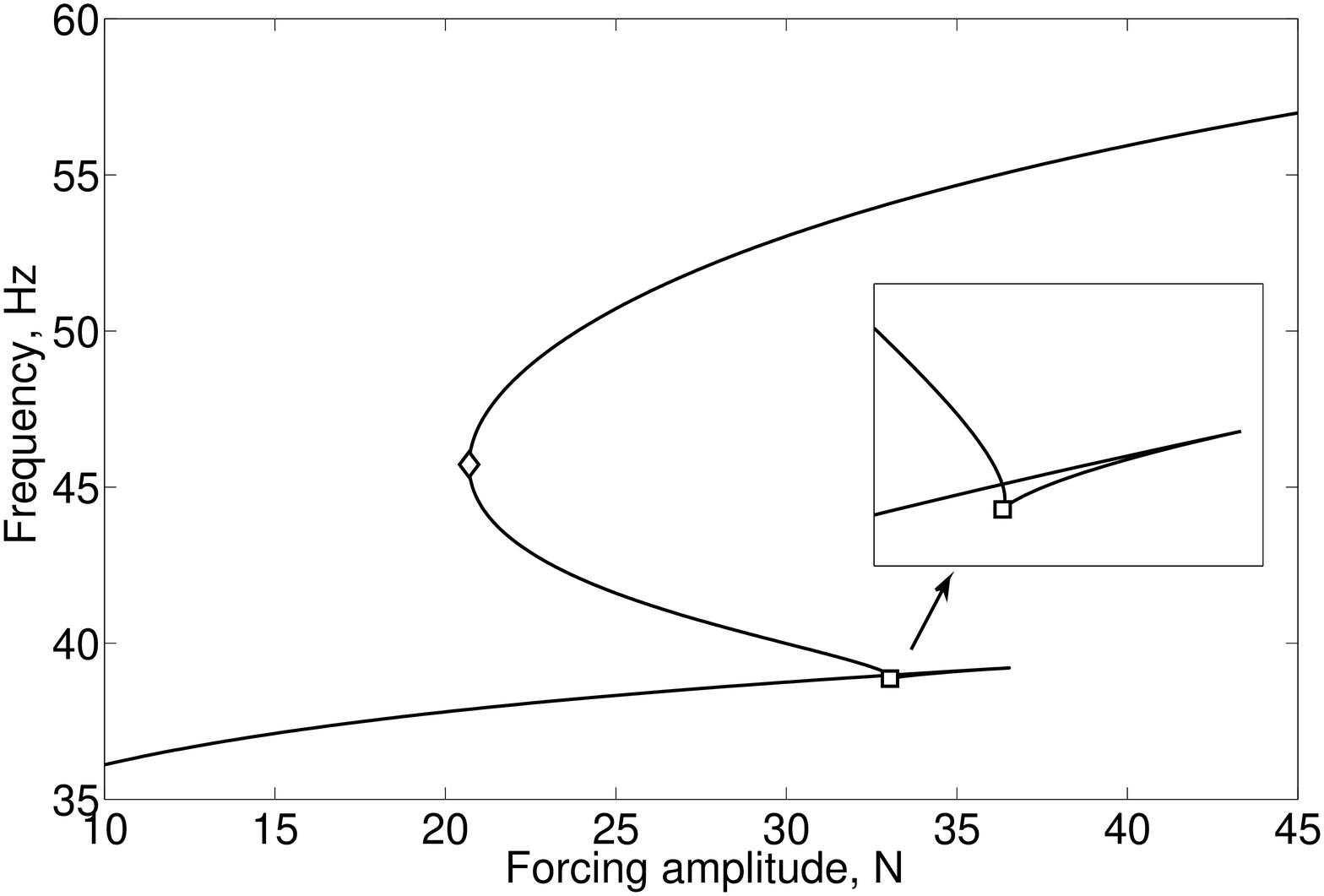}
\end{tabular}
\caption{(a) Tracking of fold bifurcations. The black lines denote nonlinear forced response curves at forcing amplitudes of 25 N, 33.1 N and 40 N, and the orange line represent the locus of fold bifurcations. (b) Projection of the bifurcation branch onto the (forcing amplitude - energy) plane. The diamond and square markers indicate the apparition and the merging of the IRC, respectively.}
\label{fig11}
\end{center}
\end{figure}

The 3D bifurcation locus is presented as an orange line in Fig. \ref{fig11}(a), which also shows the NLFR in black for forcing amplitudes of 25 N, 33.1 N and 40 N. Figure \ref{fig11}(b) gives a convenient projection of the bifurcation branch onto the frequency versus forcing amplitude plane. The turning point indicated with a diamond marker shows the frequency/forcing amplitude at which the fold bifurcations at the tips of the IRC were created. The corresponding values (20.7 N and 45.7 Hz) reflect the good predictive capability of the energy balancing criterion (22.3 N and 46.5 Hz in Section \ref{NNMsec}). The square marker indicates when the IRC merges with the main resonance peak (33.1 N).

\begin{figure}
\begin{center}
\includegraphics[width=0.75\textwidth]{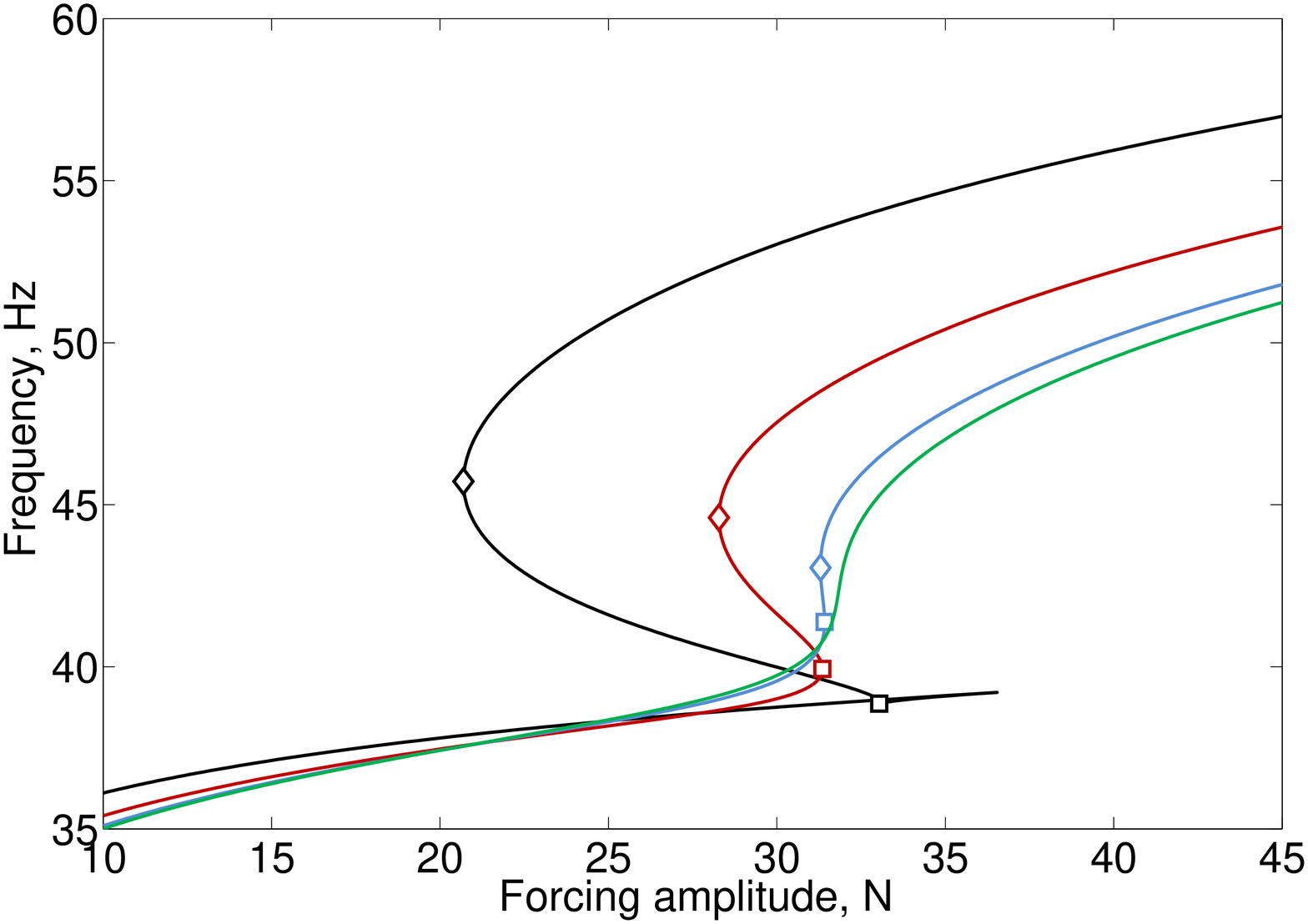}
\caption{Influence of damping on the bifurcation branches. The black, red, blue and green lines depict the branches of the system with $\kappa$ = 1, 1.5, 1.8 and 1.9, respectively. The diamond and square markers indicate the apparition and the merging of the IRC, respectively.}
\label{fig12}
\end{center}
\end{figure}

The fold bifurcation tracking analysis was also used to study the effect of structural damping on the observed IRCs. The damping matrix introduced in Section~\ref{NNMsec} was perturbed by adding a scaling term, $\kappa$, such that $\mathbf{C}=\kappa\left(-0.0391\mathbf{K}+1.47\,10^{-4}\mathbf{M}\right)$. Several bifurcation branches are given in Fig. \ref{fig12} for different values of $\kappa$, namely 1, 1.5, 1.8 and 1.9. The IRC was robust against damping since it was still visible for higher levels of damping ($\kappa>1$), however, increasing $\kappa$ caused the IRC to appear later in forcing amplitude, and shorten the range where it existed. For the largest damping case studied, for $\kappa = 1.9$, the IRC was no longer present. This analysis shows that a sufficiently large value of structural damping can destroy the IRCs. It may also explain why the other modal interactions in Fig. \ref{fig4}(a) did not produce IRCs.

\FloatBarrier

\section{Beam Response to Sine Sweep Excitation}\label{sec5}

\begin{figure}[!ht]
\begin{center}
\includegraphics[width=0.75\textwidth]{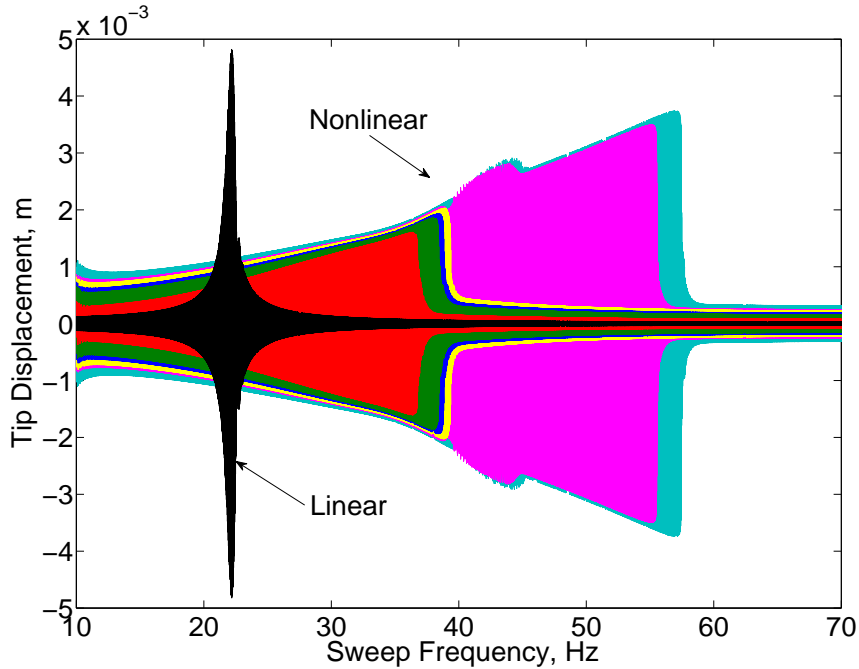}
\caption{Numerical sine sweeps at a rate of 0.5 Hz/s where the displacement of the beam tip is plotted for force amplitudes of (red) 11.3 N, (green) 22.6 N, (blue) 26.7 N, (yellow) 32.5 N, (magenta) 35.6 N and (cyan) 45.2 N.}
\label{fig5}
\end{center}
\end{figure}

To better understand the practical implications of IRCs, the beam was subjected to sine sweep excitations of increasing amplitudes with a sweep rate of 0.5 Hz/s. The outcome of the numerical simulations is plotted in Figure \ref{fig5}. For reference, the response of the linear model at a force amplitude of $A$ = 4.45 N (black) was computed, and resulted in the largest tip displacement (even though the force amplitude was lowest), with a resonance at the linear natural frequency.

For the nonlinear sweep for $A$ = 11.3 N, resonance occurred near 37 Hz, resulting in a sudden jump, the so-called jump phenomenon, to a lower response amplitude as the frequency continued to sweep upwards. When the force amplitude doubled, the resonant frequency occurred near 39 Hz. However, doubling the force amplitude once more ($A$ = 45.2 N) caused a dramatic shift in resonant frequency. Now the response dropped off around 57 Hz, indicating that the increased force amplitude shifted the resonance nearly 18 Hz. Considering the amplitudes $A$ = 32.5 N and 35.6 N shows that the shift in resonant frequency occurred in this range of forcing amplitudes. In other words, the sine sweeps at $A$ = 22.6 N, $A$ = 26.7 N and $A$ = 32.5 N fell off around 38 Hz, because the IRC was disconnected from the main branch, and there was no path for the response to follow to the higher frequency resonance. However, once the two NLFR branches merged together ($A$ = 35.6 N), the sine sweeps were able to stay along the high amplitude path up to resonance around 55 Hz.

\begin{figure}[!h]
\begin{center}
\includegraphics[width=0.75\textwidth]{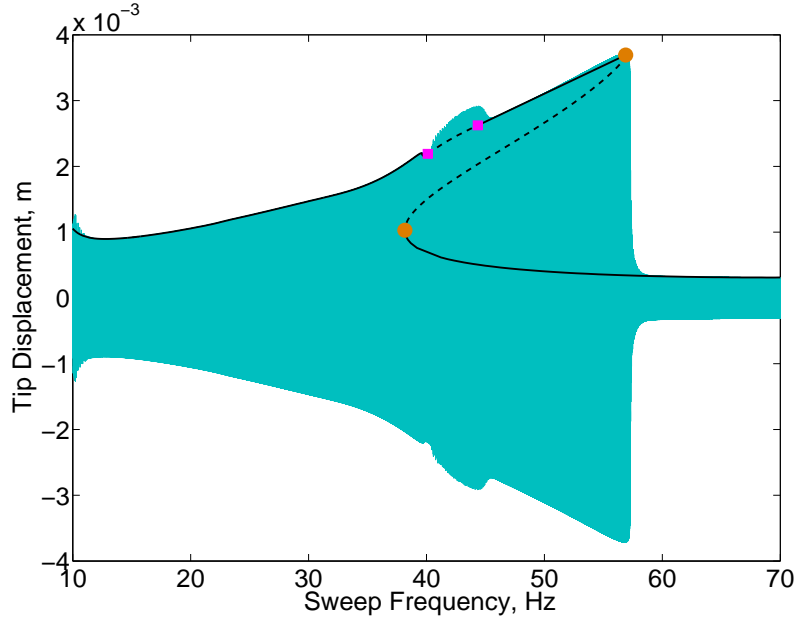}
\caption{Beam displacement at tip. Response to a swept- and a stepped-sine (numerical continuation) excitation of 45.2 N depicted in cyan and black, respectively. Solid- and dashed-black lines denote stable and unstable periodic solutions, respectively. Fold and Neimark-Sacker bifurcations are pictured with orange bullets (\textcolor{myorange}{$\bullet$}) and magenta squares (\scriptsize{\textcolor{magenta}{$\blacksquare$}}\normalsize), respectively.}
\label{fig6}
\end{center}
\end{figure}

Another dynamical phenomenon, which only appeared for $A$ = 35.6 N and 45.2 N, is the modulation of the signal's envelope in the range of 40-45 Hz. It was further examined by monitoring fold and Neimark-Sacker bifurcations~\cite{KuznetsovBook} along the NLFR. Figure~\ref{fig6} superposes the sine sweep and NLFR responses for $A$ = 45.2 N. Around 40 Hz, a Neimark-Sacker bifurcation originally located on the IRC changes the stability of the NLFR and generates a new branch of quasiperiodic oscillations (not shown in the figure). As a result, a stable torus attracts the dynamics and is responsible for the observed envelope modulation. Around 45 Hz, a second Neimark-Sacker bifurcation transforms the quasiperiodic motion back into stable periodic motion. There is a small delay between the first (second) bifurcation and the onset (disappearance) of quasiperiodic motion; this delay can be attributed to the transient character of the swept-sine excitation.

\section{Conclusion}\label{sec_future}

This paper studied the intimate connection that exists between nonlinear normal modes, i.e., the periodic motions of the undamped, unforced system, and the forced response of the damped system, with a specific focus on modal interactions and isolated resonance curves. To bridge the gap between these two types of response, the energy balancing technique was adapted to estimate the resonant response to harmonic forcing from the nonlinear modes and the damping matrix. Because it is able to reveal the presence of isolated resonance curves, the combination of nonlinear normal modes and energy balancing represents a very useful tool for global analysis of nonlinear systems.

Isolated resonance curves, which might easily be missed during numerical continuation or experimental testing, have important practical consequences for the design and testing of engineering structures. The associated response can be much larger than on the main branch, and, when they connect to the main resonance branch, they may lead to a dramatic and sudden change in resonance frequency, something which is rarely discussed in the mechanical engineering literature.

\section*{Acknowledgments}

The author RJ Kuether would like to acknowledge the funding for this part of this research from the National Physical Science Consortium (NPSC) Fellowship. The author L. Renson is Research Fellow (FRIA fellowship) of the {\it Fonds de la Recherche Scientifique - FNRS}, which is gratefully acknowledged. The authors T. Detroux, C. Grappasonni and G. Kerschen would like to acknowledge the financial support of the European Union (ERC Starting Grant NOVIB 307265).

\bibliographystyle{unsrt}


\end{document}